\documentclass[draft]{article}
\usepackage{amsmath, amssymb, latexsym}
\usepackage[T1]{fontenc}

\usepackage{amsfonts}
\usepackage{amssymb}

\newtheorem{theorem}{Theorem}

\newtheorem{corollary}[theorem]{Corollary}
\newtheorem{proposition}[theorem]{Proposition}
\newtheorem{lemma}[theorem]{Lemma}

\newtheorem{example}[theorem]{Example}

\newcommand{\der}{\operatorname{Der}}

\newcommand{\epr}{\hfill$\diamondsuit$\smallskip}
\newcommand{\pr}{\textit{Proof}. }
\title{Rings of constants of linear derivations on Fermat rings}
\author{
Marcelo  Veloso\\
e-mail:
\texttt{veloso@ufsj.edu.br}\\
Ivan Shestakov\\
e-mail:
\texttt{shestak@ime.usp.br}}
\date{}

\begin{document}
\maketitle
\begin{abstract}
\noindent In this paper we characterize all the linear $\mathbb{C}$-derivations of the Fermat ring. We show that the Fermat ring has linear $\mathbb{C}$-derivations with trivial ring of constants and construct some examples. 
\end{abstract}



\section*{Introduction}

The present paper deals with $\mathbb{C}$-derivations of the \textbf{\textit{Fermat ring}}
\[
B_n^{m} =\frac{\mathbb{C}[X_1,\ldots, X_n]}{(X_1^{m_1}+\cdots +X_n^{m_n})},
\]
where $\mathbb{C}[X_1,\ldots, X_n]$ is the polynomial ring in $n$ variables over the complex numbers $\mathbb{C}$, $n\geq3$, $m=(m_1,\ldots,m_n)$, $m_i\in \mathbb{Z}$ and $m_i \geq 2$ for $i=1,\ldots,n$.

It is well known the difficulty to describe the ring of constants of an arbitrary derivation (see \cite{BVon, FMcon, gfr, NOW}). It is also difficult to decide if the ring of constants of a derivation is trivial (see \cite{NOW, JmoAn, janusz}). In this work  we study the ring of constants of linear derivations of Fermat rings and its locally nilpotent derivations.

In $\cite{NOW}$, Andrzej Nowicki presents a description of all linear $\mathbb{C}$-derivations of the polynomial ring $\mathbb{C}[X_1,\ldots, X_n]$ which do not admit any nontrivial rational constant.

In a recent paper $\cite{BVon}$, P. Brumatti and M. Veloso show that for $m=(2,\ldots,2)$ the ring $B_n^m$ has nonzero irreducible locally nilpotent derivations. Furthermore, whenever $m_1=\cdots=m_n$, They show  that certain classes of derivations of $\mathbb{C}[X_1,\ldots, X_n]$ do not induce derivations of $B_n^m$ or are not locally nilpotent if they do.


In this work we obtain some similar results to $\cite{BVon}$, considering more general Fermat rings. We present a description of all the linear $\mathbb{C}$-derivations of $B_n^m$ when  $m=(m_1,\ldots,m_n)$ and $m_i \geq 3$ (Theorem $\ref{linerdiagonal}$)   and $m=(2,\ldots,2)$ (Theorem $\ref{da+ds}$). We also provide examples of linear derivations of $B_n^m$ with trival ring constants.

The text is organized as follows: Section 1 gathers the basic definitions and notations. Further,   we discuss several properties of the ring $B_n^m$ are discussed and a set of generators for $Der(B_n^m)$ is presented. Section 2 is dedicated to the study of the linear derivations of the Fermat ring. The set of all locally nilpotent $\mathbb{C}$-derivations of $B_n^m$ is studied in Section 3.  Finally, Sections 4 and 5 are devoted to the study of the rings of constants of linear derivations of Fermat rings.


\section{Preliminaries and Some Results}

In this paper the word \textbf{\textit{``ring"}} means a commutative ring with unit and characteristic zero. Furthermore, we denote the group of units of a ring $R$ by $R^*$ and the polynomial ring in $n$ variables over $R$ by $R[X_1,\ldots, X_n]$. A \textbf{\textit{``domain''}} is an integral domain. 


An additive mapping $D: R \rightarrow R $ is said to be a \textbf{\textit{ derivation}}  of $R$ if it satisfies the Leibniz rule: $D(ab)=aD(b)+D(a)b$, for all $a,b \in R$. If $A$ is a subring of $R$ and $D$ is a derivation of $R$ satisfying $D(A)=0$ is called $D$ an $A$-derivation. The set of all derivations of $R$ is denoted by $Der(R)$, the set of all $A$-derivations of $R$ by $Der_{A}(R)$  and by $\ker(D)$,  the \textbf{\textit{ ring of constants}} of $D$, that is $\ker(D)  =\{a \in R \mid D(r) = 0\}$. 


In this paper, the word "\textit{derivation}" \,  implicitly means a derivation which is $\mathbb{C}$-derivation and therefore we will use the notation $Der(B_n^m)$ to denote $Der_{\mathbb{C}}(B_n^m)$.  The  residue classes of variables $X$, $Y$, $Z$, ...  module an ideal are represented by  $x,$ $y$, $z$, respectively.  The symbol $\mathbb{C}$ is reserved to indicate the field of complex numbers.

A derivation $D$ is \textbf{\textit{ locally nilpotent}} if for each $r \in R$ there  is an integer $n\geq 0$ such that $D^n(r)=0$. We denote by $LND(R)$ the set of all locally nilpotent derivations of $R$. 

We say that a element $b \in R$ is a \textbf{\textit{ Darboux element}} of $D \in Der(R)$ if $b\neq 0$, $b$ is not invertible in $R$  and $D(b) = \lambda b$ for some $\lambda \in R$. In other words, a nonzero element $b$ of $R$ is a Darboux element of $D$ if, and only if, the principal ideal $(b)=\{rb \,\mid \, r\in R\}$ is different from $R$ and it is invariant with respect to $D$, that is $D((b))\subset (b)$. If $b$ is a Darboux element of $D$, then every $\lambda \in R$, such that $D(b)=\lambda b$, is said to be an \textbf{\textit{ eigenvalue}} of $b$. In particular, every element nonzero and noninvertible element belonging to the ring of constants, $\ker\, D$, is a Darboux element of $D$.   If $R$ is a domain and $D(b)=\lambda b$, then it is easy to see such the eigenvalue $\lambda$ is unique.

\begin{lemma}
\label{coef-uni}
Let $B_n^{m}$ where $m=(m_1,\,\ldots,\,m_n)$. Then, for each $f \in B_n^{m}$, there is a unique $F \in \mathbb{C}[X_1,\ldots,X_n]$ such that $deg_{X_n}<m_n$ and $f=F(x_1,\ldots,x_n)$.
\end{lemma}
\pr It follows directly from the Euclidean division algorithm by considering the polynomial $X_1^{m_1}+\cdots+X_n^{m_n}$ as a monic polynomial in $X_n$ with coefficients in $\mathbb{C}[X_1,\ldots,X_{n-1}]$.
\epr

\begin{theorem}(\cite[Theorem 4]{FMcon})
If $n\geq 5$ and $m_i \geq 2$ for all $1 \leq  i \leq  n$, then $B_n^m$ is a unique factorization domain.
\epr
\end{theorem}

We also can write $B_n^{m}=\mathbb{C}[x_1,\ldots, x_n]$, where $x_1^{m_1}+\cdots +x_n^{m_n}=0$. Here $x_1,\, x_2, \ldots,x_n$ are the images of $X_1,\,X_2,\ldots, X_n$ under the canonical epimorphism $\mathbb{C}[X_1,\ldots, X_n]\rightarrow B_n^{m}$.  An element of form $ax_1^{m_1}\cdots x_{n-1}^{m_{n-1}}$ or \linebreak $bx_1^{m_1}\cdots x_{n_1}^{m_{n-1}}x_n^j$, for  $1\leq j\leq m_n-1$, is called \textbf{\textit{monomial}}. A nonzero element $f \in \mathbb{C}[x_1,\ldots,x_n]$ is said to be \textbf{\textit{homogeneous element}} of degree $k$ if $f$ is of the form 
\[
f=\displaystyle\sum_{i_1+\cdots +i_n=k} a_{(i_1\cdots i_n)} x_1^{i_1}\cdots x_n^{i_n}
\]
where $1\leq i_n\leq m_n-1$ and $a_{(i_1\cdots i_n)}\in \mathbb{C}$ for all $(i_1\cdots i_n)$. We assume that the zero element is a homogeneous element of any degree. Furthermore, we denote by $\mathcal{V}_k$ the set of all homogeneous elements of degree $k$. Clearly $\mathcal{V}_k$ is a subspace of $B_n^m$.


\subsection{A set of generators for $Der(B_n^m)$}

Now we will present  a set of generators for the $B_n^m$-module $Der(B_n^m)$.  

First some notation will be  established. Given $H \in S = \mathbb{C}[X_1,\ldots,X_n]$ and $1\leq i \leq n$, the partial derivative $\frac{\partial (H)}{\partial X_i}$ is denoted by $H_{X_i}$. For all pairs $i,\, j \in \{1,\ldots,n\}$ with $i \neq j$, we define the derivation 
$D_{H_{ij}} = H_{X_i}\dfrac{\partial}{\partial X_j}-H_{X_j}\dfrac{\partial}{\partial X_i}$ on $S$. Observe that 
$D_{H_{ij}}(H) = 0$. 

Let $A =\frac{\mathbb{C}^{[n]}}{I}$ be a finitely generated $\mathbb{C}$-algebra. Consider the $\mathbb{C}^{[n]}$-submodule \linebreak
$
\mathcal{D}_I=\{D\in Der_\mathbb{C}(\mathbb{C}^{[n]})\mid D(I)\subseteq I\}
$
of the module $Der_\mathbb{C}(\mathbb{C}^{[n]})$. It is well known that the $\mathbb{C}^{[n]}$-homomorfism  
$\varphi : \mathcal{D}_I\rightarrow Der_\mathbb{C}(A)$ given by $\varphi(D)(g+I)=D(g)+I$ induces a  $\mathbb{C}^{[n]}$-iso\-mor\-fism of $\frac{\mathcal{D}_I}{IDer_\mathbb{C}(\mathbb{C}^{[n]})}$ in $Der_\mathbb{C}(A)$.

The Theorem $\ref{genfermat}$ will be needed, its proof can be found  in $\cite[Proposition \,\, 1 ]{BVa}$.

\begin{theorem}
\label{genfermat}
Let $F \in \mathbb{C}^{[n]} =\mathbb{C}[X_1,\ldots,X_n]$ $(n\geq 2)$ be such that $\{F_{X_1},...,F_{X_n}\}$ is a regular sequence in $S$. If there exists a derivation $\partial$ on $S$ such that $\partial(F) = \alpha F$ for some $\alpha \in \mathbb{C}$, then the $\mathbb{C}^{[n]}$-module 
\[
\mathcal{D}_F :=\{D \in Der(S) \mid D(F) \in F \cdot \mathbb{C}^{[n]}\}
\]
is generated by the derivation $\partial$ and the derivations $D_{ij} = D_{ij}^F$ for $i < j$.
\end{theorem}

From now on, the derivations $D_{F_{ij}}$, where $F = X_1^{m_1}+\cdots+X_n^{m_n}$,  will be denoted by $Dij$. Since
\[
D_{ij}(X_k)=\left\{
\begin{array}{clc}
	-m_jX_j^{m_j-1}& if & k=i\\
	m_iX_i^{m_i-1} & if &  k=j\\
	0   & if &  k\notin \{i,j\}\\
\end{array}\right.
\]
so $D_{ij}(F)=0$. Then $D_{ij} \in Der(S)$ induces $d_{ij}=m_ix_i^{m_i-1}\dfrac{\partial}{\partial x_j}-m_jx_j^{m_j-1}\dfrac{\partial}{\partial x_i}$ in $Der(B_n^m)$. Consider  the derivation
\[
E=\dfrac{1}{m_1}X_1\dfrac{\partial}{\partial X_1}+\cdots+\dfrac{1}{m_n}X_n\dfrac{\partial}{\partial X_n}.
\]
Note that $E$ satisfies $E(F) = F$. Hence, $E \in Der(S)$ induces 
\[
\varepsilon = \dfrac{1}{m_1}x_1\dfrac{\partial}{\partial x_1}+\cdots+\dfrac{1}{m_n}x_n\dfrac{\partial}{\partial x_n} \in Der(B_n^m) 
\] 

As a consequence of Theorem $\ref{genfermat}$ the following result is obtained: 

\begin{proposition}
If $F = X_1^{m_1}+\cdots+X_n^{m_n}$ then $\mathcal{D}_F :=\{D \in  Der(S)\mid D(F)\in F \cdot S\}$ is generated by the derivation $E$ and the derivations $D_{ij}$, $i <j$. In particular the $B_n^m$-module $Der(B_n^m)$ is generated by the derivation $\varepsilon$ and by the derivations $d_{ij}$, for $i < j$.
\end{proposition}
\pr
Since $\{m_1X^{m_1-1},\ldots, m_nX^{m_n-1} \}$ is a regular sequence and $E(F)=F$ the result following by Theorem $\ref{genfermat}$.
\epr


\section{Linear derivations}

This section is dedicated to the study of the linear derivations of the Fermat ring 
\[
B_n^{m}=\mathbb{C}[x_1,\ldots, x_n],
\] 
where $x_1^{m_1}+\cdots +x_n^{m_n}=0$.

A derivation $d$ of the ring $B_n^m$ is called  \textbf{linear} if 
\[
d(x_i)=\displaystyle\sum_{j=1}^n a_{ij}x_j \mbox{ for } i=1,\ldots,  n, \mbox{ where } a_{ij} \in \mathbb{C}.
\]
The matrix  $[d]=[a_{ij}]$ is called the \textbf{associated matrix } of the derivation  $d$.

\begin{theorem}
\label{linerdiagonal}  
Let $\:\:d \in Der(B_n^{m})$ be linear. If  $m=(m_1,\ldots, m_n)$ with $m_i \geq 3$ for all $i=1,\ldots, n$, then its associated matrix $[d]$ is a diagonal matrix and has the following form		
\[
\begin{bmatrix}
  \frac{\alpha}{m_1}&                   &       &\\ 
                    &\frac{\alpha}{m_2} &       &\\ 
                    &                   & \ddots&\\
	                  &                   &       & \frac{\alpha}{m_n}\\
\end{bmatrix}.
\]
for some $\alpha \in \mathbb{C}$.	
\end{theorem}
\pr
Let $[d]=[a_{ij}]$ be the associated matrix of $d$. Then $d(x_i)=\displaystyle\sum_{j=1}^n a_{ij}x_j$, for all $i$. Since $x_1^{m_1}+\cdots+x_n^{m_n}=0$,  
\[
{m_1}x_1^{m_1-1}d(x_1)+\cdots+{m_n}x_n^{m_i-1}d(x_n)=0.
\]
Then,
\begin{equation}
\label{genlinear1}
0={m_1}x_1^{m_1-1}(\sum_{j=1}^na_{1j}x_j)+{m_2}x_2^{m_2-1}(\sum_{j=1}^na_{2j}x_j)+\cdots+{m_n}x_n^{m_n-1}(\sum_{j=1}^na_{nj}x_j)
\end{equation}
Now note that  
\[
\begin{split}
m_1x_1^{m_1-1}(\sum_{j=1}^na_{1j}x_j)=&m_1a_{11}(x_1^{m_1})+m_1\sum_{j\neq 1}^na_{1j}x_jx_1^{m_1-1}\\
=&m_1a_{11}(-x_2^{m_2} -\cdots-x_n^{m_n})+m_1\sum_{j\neq 1}^na_{1j}x_jx_1^{m_1-1}\\
\end{split}
\]
and
\[
{m_2}x_2^{m_2-1}(\sum_{j=1}^na_{2j}x_j)={m_2}a_{22}x_2^{m_2}+{m_2}\sum_{j\neq 2}^na_{2j}x_jx_2^{m_2-1}
\]
\[\vdots\]
\[
{m_n}x_n^{m_n-1}(\sum_{j=1}^na_{nj}x_j)={m_n}a_{nn}x_n^{m_n}+{m_n}\sum_{j\neq n}^na_{nj}x_jx_n^{m_n-1}
\]
replacing in the Equation $\eqref{genlinear1}$ we obtain

\begin{equation}
\label{eq:genlinear}
\begin{split} 
0&=(m_2a_{22}-m_1a_{11})x_2^{m_2}+\cdots+(m_na_{nn}-m_1a_{11})x_n^{m_n}+ {m_1}\sum_{j\neq 1}^na_{1j}x_jx_1^{{m_1}-1}+\\
& {m_2}\sum_{j\neq 2}^na_{2j}x_jx_2^{{m_2}-1}+\cdots+{m_n}\sum_{j\neq n}^na_{nj}x_jx_n^{{m_n}-1}.
\end{split} 
\end{equation}

\noindent Observe that if $m_i \geq 3$, then  
\[
\{x_2^{m_1},\ldots,x_n^{m_n}\}\cup \{x_jx_i^{m_i-1}\,\,\mid\:\:1\leq i<j\leq n, \}\cup \{x_jx_i^{m_i-1}\,\,\mid\:\:1\leq j<i\leq n \}
\] 
is a linearly independent set over $\mathbb{C}$. Thus, we conclude that 
\[
m_na_{nn}=\cdots=m_2a_{22}=m_1a_{11}=\alpha \mbox{ and } a_{ij}=0 \mbox{ if } i\neq j,
\]
i.e.
\[
a_{ij}=
\left\{\begin{array}{ccc}
0& \mbox{ if } & i\neq j	\\
\frac{\alpha}{m_{i}}& \mbox{ if } & i=j
\end{array}\right.
\]
\epr

This theorem shows that for $m=(m_1,\ldots,m_n)$ and $m_i \geq 3$  linear derivations of $B_2^m$ are what is called  \textbf{\textit{ diagonal derivations}}. 

The next result characterizes linear derivations of $B_n^m$ whenever \linebreak $m=(2,\ldots, 2)$. Previously, remember that a square matrix with complex elements $A$ is said to be \textbf{\textit{skew-symmetric matrix}} if $A^T=-A$ (here $A^T$ stands, of course, for the transpose of the matrix $A$).

\begin{theorem}
\label{da+ds}
Let $\:\:d \in Der(B_n^{m})$ be linear. If  $m=(2,\ldots, 2)$, then there exist a scalar derivation $d_{\alpha}$ ($[d_{\alpha}]$ is a scalar matrix) and a skew-symmetric derivation $d_s$ ($[d_s]$  a skew-symmetric matrix) such that $d=d_{\alpha}+d_s$. This decomposition is unique.
\end{theorem}
\pr
Let  $d\in Der(B_n^m)$ be a linear derivation and $A=[a_{ij}]$ its associated matrix.
Using the same arguments used in Theorem $\ref{linerdiagonal}$ we obtain  

\[
0=(a_{22}-a_{11})x_2^2+\cdots+(a_{nn}-a_{11})x_n^2 + \sum_{i<j}(a_{ij}+a_{ji})x_ix_j
\]
Since the set $\{x_2^2,\ldots,x_n^2\}\cup \{x_ix_j;\:\:1\leq i<j\leq n\}$ is linearly independent over $\mathbb{C}$, it follows that  
\[
a_{11}=a_{22}=\cdots =a_{nn}=\alpha \:\:\mbox{and}\:\: a_{ij}=-a_{ji} \:\:\mbox{if}\:\: i<j,
\]
then its associated matrix $[d]$ has the following form		
\[
\begin{bmatrix}
            \alpha  & a_{12}  & \ldots  & a_{1n}\\ 
        -a_{12}&  \alpha      &         & a_{2n}\\ 
         \vdots&  \vdots & \ddots  & \vdots\\
	      -a_{1n}& -a_{2n} & \ldots  & \alpha
\end{bmatrix}.
\]
where $\alpha, a_{ij} \in \mathbb{C}$. Now define $d_{\alpha}$ by $d_{\alpha}(x_i)=\alpha x_i$, $i=1,\ldots, n$ and $d_s=d-d_{\alpha}$.

\epr


\section{Locally Nilpotent Derivations}

In this section we proof that the unique locally nilpotent derivation linear of $B_n^m$ for $m=(m_1,\ldots,m_n)$ and $m_i \geq 3$ is the zero derivation. Further,  we show that a certain class of derivations of $\mathbb{C}[X_1,\ldots, X_n]$ do not induce nonzero locally nilpotent derivation  of $B_n^m$.

Let $S =\frac{\mathbb{C}^{[n]}}{I}$ be a finitely generated $\mathbb{C}$-algebra. Consider the $\mathbb{C}^{[n]}$-submodule \linebreak
$
\mathcal{D}_I=\{D\in Der_\mathbb{C}(\mathbb{C}^{[n]})\mid D(I)\subseteq I\}
$
of the module $Der_\mathbb{C}(\mathbb{C}^{[n]})$. It is well known that the $\mathbb{C}^{[n]}$-homomorfism  
$\varphi : \mathcal{D}_I\rightarrow Der_\mathbb{C}(S)$ given by $\varphi(D)(g+I)=D(g)+I$ induces a \linebreak  $\mathbb{C}^{[n]}$-iso\-mor\-fism of $\frac{\mathcal{D}_I}{IDer_\mathbb{C}(\mathbb{C}^{[n]})}$ in $Der_\mathbb{C}(S)$. From this fact  the following result is obtained.

\begin{proposition}
\label{gentriangular}
Let $d$ be a derivation of the $B_n^{m}$. If $d(x_1)=a\in \mathbb{C}$ and for each $i$, $1<i\leq n$, $d(x_i)\in \mathbb{C}[x_1,\ldots,x_{i-1}]$, then $d$ is the zero derivation.
\end{proposition}

\begin{pr} 
Let $F$ be the polynomial $X_1^{m_1}+\cdots +X_n^{m_n}$. We know that exists $D\in Der(\mathbb{C}^{[n]})$ such that $D(F)\in F\mathbb{C}^{[n]}$ and that  $d(x_i)=D(X_i)+ F\mathbb{C}^{[n]}$, $\forall i$. Thus  $D(X_1)-a\in F\mathbb{C}^{[n]}$, and for each  $i>1$ there exists  $G_i=G_i(X_1,\ldots,X_{i-1})\in \mathbb{C}[X_1,\ldots,X_{i-1}]$ such that $D(X_i)-G_i\in F\mathbb{C}^{[n]}$. Since $D(F)=\displaystyle\sum_{i=1}^{n}m_iX_i^{m_i-1}D(X_i)\in F\mathbb{C}^{[n]}$ and  
\[
D(F)=\displaystyle\sum_{i=1}^{n}m_iX_i^{m_i-1}(D(X_i)-G_i)+\displaystyle\sum_{i=1}^{n}m_iX_i^{m_i-1}G_i,
\]
where $G_1=a$, we obtain $\displaystyle\sum_{i=1}^{n}m_iX_i^{m_i-1}G_i\in F\mathbb{C}^{[n]}$ and then obviously $G_i=0$  for all $i$. Thus $d$ is the zero derivation.
\end{pr}

\begin{lemma}\label{limanil}
Let $d$ be a linear derivation of $B_n^m$ and $[a_{ij}]$ its associated matrix. Then $d$ is locally  nilpotent if and only if $[a_{ij}]$ is nilpotent.
\end{lemma}
\pr
The following equality can be verified by induction over $s$.  
 
\begin{eqnarray}\label{dsiax}
\left[
\begin{array}{c}
d^s(x_1)\\
\vdots \\
d^s(x_n)
\end{array}
\right]
&=&
{[a_{ij}]}^s
\left[
\begin{array}{c}
x_1 \\
\vdots \\
x_n
\end{array}
\right].
\end{eqnarray}
We know that $d$ is locally nilpotent if and only if there exists $r\in \mathbb{N}$ such that $d^r(x_i)=0$ for all $i$. As $\{x_1,\ldots,x_n\}$ is linearly independent over $\mathbb{C}$ by the \eqref{dsiax}, the result follows.
\epr

\begin{theorem}
If  $\:\:d \in LND(B_n^{m})$ is linear and $m=(m_1,\ldots, m_n)$ wich $m_i \geq 3$, then $d$ is the zero derivation.				
\end{theorem}
\pr
Since   $d$ is locally nilpotent,  $[d]$ is nilpotent (by Lemma $\ref{limanil}$) and diagonal (by Theorem $\ref{linerdiagonal}$). Thus, the matrix $[d]$ is null and  $d$ is the zero derivation.
 
\epr

In the case  $m=(2,\ldots,2)$, linear locally nilpotent derivations of the ring $B_n^m$ were characterized as follows. 

\begin{theorem}\cite[Theorem 1]{BVon}
\label{lndiffnil}
If $d \in Der(B_n^m)$ is linear and $m=(2,\ldots,2)$, then $d\in LND(B_n^2)$ if, and only if, its associated matrix is  nilpotent and skew-symmetric.\epr
\end{theorem}


\section{Ring of constants}

In this section we show that the ring of constants of  all nonzero linear derivations of $B_n^m$, where  $m=(m_1,\ldots,m_n)$ and $m_i \geq 3$, is trivial, that is $\ker(d)=\mathbb{C}$.

During all this section we always consider  $m=(m_1,\ldots,m_n)$ with $m_i \geq 3$. 

The next result ensures the existence of  Darboux elements for every nonzero linear derivation of  $B_n^{m}$.

\begin{proposition}
\label{darbouxmon}
Let $d$ be a nonzero linear derivation of $B_n^{m}$. If $d(x_i)=\frac{\alpha}{m_i}x_i$, $i=1,\ldots, n$, for some $\alpha \in \mathbb{C}^*$, then  $f=b x_1^{i_1}\cdots x_n^{i_n}$ is a Darboux element of $d$, $d(f)=\lambda f$, and 
\[
\lambda=\alpha \left(\frac{i_1}{m_1}+\frac{i_2}{m_2}+\cdots+\frac{i_n}{m_n}\right).
\]
\end{proposition}
\pr
Let $f= b x_1^{i_1}\cdots x_n^{i_n}$. Then
\[
\begin{split}
d(f)
=&d(b x_1^{i_1}\cdots x_n^{i_n})\\
=&b d(x_1^{i_1}\cdots x_n^{i_n})\\
=&b\displaystyle\sum_{k=1}^n  i_kx_1^{i_1}\cdots x_k^{i_k-1}\cdots x_n^{i_n}d(x_k)\\
=&b\displaystyle\sum_{k=1}^n  i_kx_1^{i_1}\cdots x_k^{i_k-1}\cdots x_n^{i_n}\left(\frac{\alpha}{m_k}x_k\right)\\
=&\alpha b\displaystyle\sum_{k=1}^n  \frac{i_k}{m_k}x_1^{i_1}\cdots  x_n^{i_n}\\
=&b\alpha \left(\frac{i_1}{m_1}+\frac{i_2}{m_2}+\cdots+\frac{i_n}{m_n}\right) x_1^{i_1}\cdots x_n^{i_n}\\
=&\lambda f\\
\end{split}
\]
\epr

\begin{corollary}
Let $d$ be a nonzero linear derivation of $B_n^{m}$. If $f$ is a homogeneous element of degree $k$, then $f$ is a Darboux element of  $B_n^m$ with eigenvalue $\lambda =\frac{k}{m}$.
\end{corollary}
\pr
Let $f=\displaystyle\sum_{i_1+\cdots +i_n=k} a_{(i_1\cdots i_n)} x_1^{i_1}\cdots x_n^{i_n}$ be a homogeneous element of degree $k$, where $0\leq i_n <m$ and $a_{(i_1\cdots i_n)}\in \mathbb{C}$.

\[
\begin{split}
d(f)
=&d\left(\displaystyle\sum_{i_1+\cdots +i_n=k} a_{(i_1\cdots i_n)} x_1^{i_1}\cdots x_n^{i_n}\right)\\
=&\displaystyle\sum_{i_1+\cdots +i_n=k} a_{(i_1\cdots i_n)} d(x_1^{i_1}\cdots x_n^{i_n})\\
=&\displaystyle\sum_{i_1+\cdots +i_n=k} a_{(i_1\cdots i_n)} \left(\frac{i_1}{m}+\frac{i_2}{m}+\cdots
   +\frac{i_n}{m}\right)x_1^{i_1}\cdots x_n^{i_n}\\
=&\displaystyle\sum_{i_1+\cdots +i_n=k} a_{(i_1\cdots i_n)} \frac{k}{m} x_1^{i_1}\cdots x_n^{i_n}\\
=&\frac{k}{m}\left(\displaystyle\sum_{i_1+\cdots +i_n=k} a_{(i_1\cdots i_n)} x_1^{i_1}\cdots x_n^{i_n}\right)\\
=&\lambda f\\
\end{split}
\]

\epr

The main result this section is:

\begin{theorem}
\label{kerd=c}
Let $d$ be a nonzero linear derivation of $B_n^{m}$. Then $\ker(d)=\mathbb{C}$.
\end{theorem}
\pr
By Theorem $\ref{linerdiagonal}$, $d(x_i)=\frac{\alpha}{m_i}x_i$ for  $i=1,\ldots, n$ and $\alpha \in \mathbb{C}$. Since $d \neq 0$, so $\alpha \neq 0$. Let $0 \neq f \in B_n^m$ such that $d(f)=0$. Thus $f=\displaystyle\sum_{(i_1,\ldots,i_n) \in I} b_{(i_1,\ldots,i_n)}x_1^{i_1}\cdots x_n^{i_n}$ where  $0\neq b_{(i_1,\ldots,i_n)} \in \mathbb{C}$ for all $(i_1,\ldots,i_n) \in I$. Then
\[
\begin{split}
0=d(f)
=&\displaystyle\sum b_{(i_1\cdots i_n)}d(x_1^{i_1}\cdots x_n^{i_n})\\
=&\displaystyle\sum b_{(i_1\cdots i_n)}\alpha\left(\frac{i_1}{m_1}+\frac{i_2}{m_2}+\cdots+\frac{i_n}{m_n}\right)x_1^{i_1}\cdots x_n^{i_n}\\
\end{split}
\]
It follows from Lemma $\ref{coef-uni}$ that $b_{(i_1\cdots i_n)}\alpha(\frac{i_1}{m_1}+\frac{i_2}{m_2}+\cdots+\frac{i_n}{m_n} )\neq 0$   for all $(i_1,\,\ldots, i_n)\in I$, because $b_{(i_1\cdots i_n)}\alpha \neq 0$ for all $(i_1,\,\ldots, i_n)\in I$. So  $\frac{i_1}{m_1}+\frac{i_2}{m_2}+\cdots+\frac{i_n}{m_n} =0$ for all $(i_1,\,\ldots, i_n)\in I$. This implies that $(i_1,\ldots, i_n)=(0,\ldots,0)$ for all $(i_1,\ldots,i_n) \in I$. Therefore $f \in \mathbb{C}$.
\epr

\begin{theorem}
Let $d \in \der(B_n^m)$ be given by $d(x_i)=\frac{\alpha}{m_i}$, $i=1,\ldots, n$, for some $\alpha \in \mathbb{C}$. If 
\[
f=\displaystyle\sum_{(i_1,\ldots,i_n) \in I} a_{(i_1,\ldots,i_n)}x_1^{i_1}\cdots x_n^{i_n}
\]
is  a darboux element of $d$, this is,  $d(f)=\lambda f$ for some $\lambda \in B_n^m$, then \linebreak $\lambda=\alpha(\frac{i_1}{m_1}+\frac{i_2}{m_2}+\cdots+\frac{i_n}{m_n} )$ for all $(i_1,\ldots,i_n) \in I$. 
\end{theorem}

\pr
Let $f=\displaystyle\sum a_{(i_1,\ldots,i_n)}x_1^{i_1}\cdots x_n^{i_n}$. 
It follows from  Theorems $\ref{kerd=c}$ and $\ref{linerdiagonal}$  that 
\[
d(f)=\displaystyle\sum_{(i_1,\ldots,i_n) \in I} a_{(i_1\cdots i_n)}b_{(i_1 \ldots i_n)} x_1^{i_1}\cdots x_n^{i_n}
\]
where $b_{(i_1\cdots i_n)}=\alpha(\frac{i_1}{m_1}+\frac{i_2}{m_2}+\cdots+\frac{i_n}{m_n} )\in \mathbb{C}$.
Then  
\[
 \displaystyle\sum a_{(i_1\cdots i_n)}b_{(i_1 \ldots i_n)} x_1^{i_1}\cdots x_n^{i_n}
=d(f)=\lambda f=\displaystyle\sum \lambda a_{(i_1\cdots i_n)}x_1^{i_1}\cdots x_n^{i_n}
\]
So
\[
a_{(i_1 \ldots i_n)}b_{(i_1\cdots i_n)}= \lambda a_{(i_1 \ldots i_n)}
\]
for all $(i_1,\ldots,i_n) \in I$, by Lemma $\ref{coef-uni}$. Therefore  $\lambda=b_{(i_1 \ldots i_n)}=\alpha(\frac{i_1}{m_1}+\frac{i_2}{m_2}+\cdots+\frac{i_n}{m_n} )$ for all $(i_1,\ldots,i_n) \in I$. 
\epr


\section{The case $m=(2,\ldots,2)$}

In this section we focus on the Fermat rings
\[
B_n^m=\frac{\mathbb{C}[X_1,\ldots, X_n]}{(X_1^{2}+\cdots +X_n^{2})},
\]
where $m= (2,\ldots,2)$ and $n\geq 3$. For simplicity we denote $B_n^m$ by $B_n^2$.

We study  linear derivations of $B_n^2$, their  rings of constants, and we show how to construct examples of  linear derivations with trivial ring of constants.

The next result  will be useful for this purpose. 

\begin{proposition}
\label{df=lamkaf}
Let $d$ be a  nonzero linear derivation of $B_n^m$, with $d=d_{\alpha}+d_s$, where  $d_{\alpha}$ is the scalar derivation and is skew-symmetric derivation. Let $d_{\alpha}$ definid by $d_{\alpha}(x_i)=\alpha x_i$ for $i=1,\ldots, n$ and  $\alpha \in \mathbb{C}$. If  $f \in B_n^2$ is homogeneous element of degree $k$ then $d(f)=\lambda f$ if, and only  if, $d_s(f)=(\lambda-k\alpha)f$.
 \end{proposition}
\pr 
It is easy to see that $d_{\alpha}(f)=kaf$. Suppose that $d(f)=\lambda f$. Then 
\[
\lambda f =d(f)=d_{\alpha}(f)+d_s(f)=-\alpha kf+d_s(f).
\]
Hence,
\[
d_s(f)=(\lambda -k\alpha )f.
\] 
Now suppose $d_1(f)=(\lambda -ka)f$. Then 
\[
d(f)=d_{\alpha}(f)+d_s(f)=k\alpha f+(\lambda -k\alpha )f= \lambda f.
\]
\epr

\begin{corollary}
\label{ds=mka}
Let $d=d_{\alpha}+d_s$ be a  nonzero linear derivation of $B_n^2$, where  $d_{\alpha}$ is the scalar derivation and $d_s$ is skew-symmetric derivation.  If  $f \in B_n^2$ is homogeneous element of degree $k$ then $d(f)=0$ if, and only if, $d_s(f)=-k \alpha f$.
 \end{corollary}
\pr 
Consider $\lambda=0$ in  Proposition $\ref{df=lamkaf}$.
\epr

The next Theorem shows that every skew-symmetric derivation has a nontrivial ring of  constants and every nonzero scalar derivation has trivial ring of constants.  

\begin{theorem}
Let $d=d_{\alpha}+d_s$ be a  nonzero linear derivation of $B_n^2$, where  $d_{\alpha}$ is the scalar derivation and $d_s$ is skew-symmetric derivation. Then

\begin{enumerate}
	\item If $d_{s}$ is zero the derivation, then $\ker(d)=\ker(d_{\alpha})$ is trivial.

  \item If $d_{\alpha}$ is the zero  derivation, then $\ker(d)=\ker(d_s)$ is nontrivial.
\end{enumerate}
 \end{theorem}
\pr 1) Observe that $d_{\alpha}(f)=kaf$ for all homogeneous element of degree $k$ and $d_{\alpha}({\mathcal{V}_k})\subset \mathcal{V}_k$. 

2) It suffices to prove that there is $f\in B_n^2$ such that $d_s(f)=0$ and $f \not\in \mathbb{C}$. Let $f$ a homogeneous element of degree $2$ of $B_n^2$, then $f=XAX^T$ where $A=[a_{ij}]$ is a symmetric matrix and $X=(x_1,\ldots,x_n)$. Observe that for  $B=[d_s]$ we have  
\[
\begin{split}
d(f)=d(XAX^T)
=&(XB)AX^T+XA(XB)^T\\
=& XBAX^T+XA(-BX^T)\\
=& XBAX^T-XABX^T\\
=& X(BA-AB)X^T.\\
\end{split}
\]
If $B^2\neq 0$ consider the symmetric matrix  $B^2$. It follows from the above remark that for $f=XB^2X^T$  we have $d_s(f)=0$ and $f \not\in \mathbb{C}$, because $A=B^2\neq 0$. If $B^2=0$, then $\lambda=0$ is eigenvalue of $B$ and $B^T$. In this case, choose a nonzero element $f=a_1x_1+\cdots + a_nx_n$ such that the nonzero vector $(a_1,\ldots, a_n)^T$ is an eigenvector of $B^T$. So $d(f)=0$ and $f \not\in \mathbb{C}$. Therefore, $\ker(d_s)\neq \mathbb{C}$.
\epr

We also show that there are linear derivations with $d_{\alpha} \neq 0$, $d_s\neq 0$, and trivial ring of constants. 

\begin{theorem}
\label{ker=c}
Let $d_s$ be a  nonzero skew-symmetric derivation of $B_n^2$.  Then exists a scalar derivation $d_{\alpha}$ of $B_n^2$ such that the derivation $d=d_{\alpha}+d_s$ satisfies $\ker(d)=\mathbb{C}$.
 \end{theorem}
\pr 
First note that the vector space $\mathcal{V}_k$ (the set of homogeneous elements  of degree $k$ of $B_n^2$) is invariant with respect to $d_s$. Hence $d_s(f)=0$ if only if $d_s(f_k)=0$ for all homogeneous components $f_k$ of $f$. As a consequence of this fact we assume that $f$ is a homogeneous element of degree $k$. Let  $\alpha$ be a nonzero complex number that satisfies the conditions
\begin{enumerate}
	\item $\alpha \notin Spec(d_s)$,
	\item for all positive integer $k$, $-k\alpha \notin Spec(d_s\mid \mathcal{V}_k)$.
\end{enumerate}
This number exists because $\mathbb{C}$ is uncountable and the set of the numbers that  satisfies $1)$ and $2)$ are countable. Let $\alpha$ be a number that satisfies the conditions $1)$ and $2)$, then $d_s(f)\neq \alpha f$ and $d_s(f)\neq -k\alpha f$,  for all $f \notin \mathbb{C}$ and for all  positive integer $k$ . Let $d_{\alpha}$ be a scalar derivation defined by $d_{\alpha}(x_i)=\alpha x_i$, $i=1,\ldots,n$. 

Finally, by considering the derivation $d=d_{\alpha}+d_s$, we show that  $\ker(d)=\mathbb{C}$. In order to do that, if $g \in B_n^2$ is a nonzero homogeneous element of degree $k$, then $d(g)=0$ if, and only if, $d_s(g)=-k\alpha g$, by Corollary $\ref{ds=mka}$, which implies $k=0$. Therefore, $g \in \mathbb{C}^*$.
\epr

We now provide and explicit example of such derivation:

\begin{example}
\label{exnonlnd}
Let $d=d_1+d_s$ be the linear derivation of $B_3^2=\mathbb{C}[x,y,z]$ given by 
\[
[d]=[d_1]+[d_s]
=\left[
\begin{array}{crr}
 1 & 0 & 0 \\
 0 & 1 &-1 \\
 0 & 1 & 1 \\
 \end{array}
\right]
=
\left[
\begin{array}{crr}
 1 & 0 & 0 \\
 0 & 1 & 0 \\
 0 & 0 & 1 \\
 \end{array}
\right]
+
\left[
\begin{array}{crr}
 0 & 0 & 0 \\
 0 & 0 &-1 \\
 0 & 1 & 0 \\
 \end{array}
\right].
\] 
We claim that  $\ker(d) = \mathbb{C}$. To be more precise, let $\mathcal{V}_k$ be the set of all homogeneous elements of degree $k$ in $\mathbb{C}[y,z]$. The set  
\[
S_k=\{y^k, y^{k-1}z,\ldots, yz^{k-1}, z^k \}
\]
is a basis for $\mathcal{V}_k$. The matrix
\[
[d\mid \mathcal{V}_k]=\left[
\begin{array}{rccrcccr}
 k       & 1      & 0      & \ldots&        &\ldots& 0     & 0 \\
 -k      & k      & 2      & \ldots&        &\ldots& 0     & 0\\
 0       & -(k-1) & k      & \ddots&        &      &\vdots & \vdots \\
         &  0     &-(k-2)  & \ddots&        & k-2  &   0   & 0 \\
 \vdots  & \vdots &  0     & \ddots&        &  k   & k-1   & 0\\
0        & 0      & \vdots &       &        &  -2  &    k  & k\\
 0       & 0      & 0      & \ldots&        & 0    &   -1  & k\\
 \end{array}
\right]
\]
is the matrix the linear derivation $d$ restrict to subspace $\mathcal{V}_k$ in the basis $S_k$. It is easy check that $Det([d\mid \mathcal{V}_k]) \neq 0$ for all  $k\geq 1$, by the principle of induction. Then $d(f)\neq 0$ for all  homogeneous elements of degree $k\geq 1$. Therefore,   $\ker(d)=\mathbb{C}$.
\end{example}

\begin{theorem}
\label{dslnd}
Let $d=d_{\alpha}+d_s$ be a  nonzero linear derivation of $B_n^2$, where  $d_{\alpha}$ is a scalar derivation and $d_s$ is a skew-symmetric derivation. If $d_s$ is a locally nilpotent derivation then  $\ker(d) =\mathbb{C}$, for all  nonzero scalar derivation $d_{\alpha}$.
 \end{theorem}
\pr 
Let $0\neq f \in B_n^2$ such that $d(f)=0$. It suffices  to show that $f \in \mathbb{C}$. We may assume that $f$ is a nonzero homogeneous element of degree $k$, because $V_k$ is invariant by $d$. Let $m$ be the smallest positive integer such that $g=d_s^{m-1}(f)\neq 0$ and $d_s^m(f)=0$, this $m$  exists because $d_s$ is a locally nilpotent derivation. This implies that $g$ is a nonzero homogeneous element of degree $k$, because $V_k$ is invariant by $d$. One easily verifies that $d_{\alpha}d_s=d_sd_{\alpha}$ and $dd_s=d_sd$. Note that $d(g)=d_{\alpha}(g)+d_s(g)=k\alpha g$, because $g$ is an homogeneous element of degree $k$ and $d_s^m(f)=0$. Now observe that
\[
d(g)=d(d_s^{m-1}(f))=d_s^{m-1}(d(f))=d_s^{m-1}(0)=0.
\]
We thus get $k\alpha g=0$. And so $k=0$. Therefore, $f \in \mathbb{C}$.
\epr

To conclude this section  we construct two families of examples which illustrate this theorem.

\begin{example}
\label{eximpar}
Let $n\geq 3$ be an odd number and $d_s$ a linear derivation of  $B_n^2$ defined by the skew-symmetric matrix $n\times n$
\[
[d_s]=
\left[
\begin{array}{cccccc}
 0&0&\ldots&0&0&-1\\
 0&0&\ldots&0&0&-i\\
 \vdots&\vdots&\ddots&\vdots&\vdots&\vdots\\
 0&0&\ldots&0&0&-1\\
 0&0&\ldots&0&0&-i\\
 1&i&\ldots&1&i&0
\end{array}
\right].
\]
It is easy to check that $[d_s]^3=0$, which implies that $[d_s]$ is nilpotent. Then $d_s$ is a locally nilpotent linear derivation of $B_n^2$,  by Theorem $\ref{lndiffnil}$. Now consider the linear derivation $d=d_1+d_s$, where $d_1(x_i)=x_i$ for $i=1,\ldots, n$. It follows from Theorem $\ref{dslnd}$ that $\ker(d_s)=\mathbb{C}$.
\end{example}

\begin{example}\label{expar}
Let $n\geq 4$ be an even number and $\varepsilon \in \mathbb{C}$ a primitive \linebreak $(n-1)$-th root of unity. Set $d_s$ a linear derivation of  $B_n^2$  by the skew-symmetric matrix $n\times n$
\[
[d_s]=\left[
\begin{array}{ccccccc}
 0 & 0 & \ldots&0&\ldots & 0 & -1 \\
 0 & 0 & 0 &\ldots& 0 & 0 & -\varepsilon \\
 \vdots &\vdots&\ddots&\vdots&\ddots&\vdots&\vdots \\
 0 & 0 & \ldots & 0 & \ldots &0 & -{\varepsilon}^k \\
 \vdots &\vdots&\ddots&\vdots&\ddots&\vdots&\vdots \\
 0 & 0 & \ldots &0 & \ldots &0& -{\varepsilon}^{n-2} \\
 1& \varepsilon & \ldots &{\varepsilon}^{k}&\ldots &{\varepsilon}^{n-2}&0
\end{array}
\right]
\]
Again, $[d_s]$ is nilpotent ($[d_e]^3=0$). Thus, $d_s$ is a locally nilpotent derivation of $B_n^2$, by Theorem $\ref{lndiffnil}$. Now considering the linear derivation $d=d_1+d_e$, where $d_1$ is the same as in the previous example, we conclude that $\ker(d_s)=\mathbb{C}$, by Theorem $\ref{dslnd}$.
\end{example}

\textbf{Remark}:
In the  Example $\ref{exnonlnd}$ it easy see that $[d_s]$ is not nilpotent and, consequently  $d_s$ is not locally nilpotent (Theorem $\ref{lndiffnil}$). This shows that $d_s$ locally nilpotent is not a necessary condition in Theorem $\ref{dslnd}$ for $\ker(d_{\alpha}+d_s)=\mathbb{C}$.


\end{document}